\documentclass{amsart}

\newcommand{\remark}{\noindent {\bf Remark. }}

\newcommand{\ws}{\hspace{4pt}}
\newtheorem{theorem}{Theorem}

\newtheorem{lemma}{Lemma}
\newtheorem{defi}{Definition}
\usepackage[]{amssymb,amsopn}
\usepackage[]{amsmath}
\usepackage[]{eucal}
\usepackage[]{latexsym}

\begin{document}

\title[Exceptional Laguerre and Jacobi Polynomials]{The Electrostatic Properties of Zeros of Exceptional Laguerre and Jacobi Polynomials and stable interpolation}
\author{\'A. P. Horv\'ath }

\thanks{Supported by Hungarian National Foundation for Scientific Research, Grant No.  K-100461}

\begin{abstract}   
We will examine the electrostatic properties of exceptional and regular zeros  of $X_m$-Laguerre and $X_m$-Jacobi polynomials. Since there is a close connection between the electrostatic properties of the zeros and the stability of interpolation on the system of zeros, we can deduce an Egerv\'ary-Tur\' an type result as well. The limit of the energy on the regular zeros are also investigated.

\end{abstract}
\maketitle

\section{Introduction}

Classical orthogonal polynomials can be introduced as the eigenfunctions of Sturm-Liouville operators and they also play a fundamental role in the construction of bound-state solutions to exactly solvable potentials in quantum mechanics. Moreover an equilibrium problem for the logarithmic interaction of positive unit charges under an external field leads to a nice electrostatic interpretation of their zeros. The following approach of generalization was investigated in the '90-s: let $w=e^{-Q}$ be a positive weight function supported on $(a,b)\subset \mathbb{R}$ with finite moments and $p_n=p_n(w)$ are the orthogonal polynomials on $(a,b)$ with respect to $w$. Assuming $Q$ is a twice differentiable and convex function, $p_n$ satisfies the following differential equation on $(a,b)$:
$$p_n^{''}(x)+M_n(x)p_n^{'}(x)+N_n(x)p_n(x)=0,$$
where $M_n(x)$ depends on $n$, actually on $p_n$ (cf. e.g. \cite{m} Th. 3.4. and \cite{iw}). Besides the classical properties of orthogonal polynomials, the electrostatic behavior of their zeros were also investigated (cf. \cite{i}). On the other hand these type of orthogonal polynomials are not the suitable ones for constructing solvable potentials. Recently some new families of orthogonal polynomials are investigated which are very useful to this purpose. At first $X_1$-Jacobi and $X_1$-Laguerre polynomials, as exceptional orthogonal polynomial families were introduced by D. G\'omez-Ullate, N. Kamran and R. Milson (cf. e.g. \cite{gukm2}). The relationship between exceptional orthogonal polynomials and the Darboux transform is observed by C. Quesne (cf. e.g. \cite{q}). Higher-codimensional families were introduced by S. Odake and R. Sasaki \cite{os}. The location of zeros of exceptional orthogonal polynomials are described by D. G\'omez-Ullate, F. Marcell\'an and R. Milson \cite{gumm}, and the electrostatic interpretation of zeros of $X_1$-Jacobi polynomials is given by D. Dimitrov and Yen Chi Lun \cite{dy}. 

Below we will extend some results of \cite{dy} to $X_m$-Laguerre polynomials of the first kind, and we will show that the regular zeros of $X_m$-Jacobi and $X_m$-Laguerre polynomials behaves like the zeros of the classical orthogonal polynomials in point of energy. To this purpose we adapt some methods was developed to general orthogonal polynomials, to the exceptional ones. Since the regular zeros of exceptional polynomials form a minimal energy (or Fekete) system with respect to a suitable external field, similarly to \cite{ho} we will show that on these sets one can build up stable interpolatory operators which are the most economical as well. Finally the notion of Fekete sets and $n^{th}$ transfinite diameter allows to investigate the behavior of the energy function when the number of the points tends to infinity.

\medskip

\section{Notations, Preliminary Computations}

We will examine an energy problem with an external field on a finite or infinite interval of $\mathbb{R}$. To this end we introduce some notations.

Let $U_n:=\{u_1, \dots , u_n\}$ be any system of nodes on an interval $I$ and $0\leq w \in C^2(I)$ be a weight function on $I$. Let $\omega_{U_n}(x):=\prod_{k=1}^n(x-u_k)$. 

\begin{defi} The energy function on $I$ with respect to $w$ is 
$$T_w(u_1, \dots , u_{n})=\prod_{j=1}^{n}w(u_j) \prod_{1\leq i<j\leq n}(u_i-u_j)^2,$$ 
cf. e.g. \cite{st} p. $143.$ or \cite{i} $(2.5)$.\end{defi}

\medskip

\remark

\noindent (1) For $w\equiv 1$ the energy function was investigated by I. Schur \cite{sch}. He called it as "discriminant" and he found that the maximum with respect to the nodes is attained at the zeros of certain (orthogonal) polynomials.

\noindent (2) L. Fej\'er recognized that the solutions of the maximum problem are adequate systems of nodes for interpolation (normal and $\varrho$-normal point systems). For the weighted version cf. \cite{aho}.

\noindent (3) It's potential theoretic meaning allows us to call $T_w$ as energy function. Indeed $\inf_{U_n\subset H}-\log \left(T_w^{\frac{1}{n(n-1)}}(u_1, \dots , u_{n})\right)$ is the discrete minimal energy or the $n^{th}$ transfinite diameter of $H$ with respect to the weight $w^{\frac{1}{2(n-1)}}$ (cf. \cite{st}).

\noindent (4) If $w$ is a classical (Jacobi, Laguerre, Hermite) weight, then the solutions of the weighted energy problem are the $1(w)$-normal systems which are the zeros of the classical orthogonal polynomials, namely for $w=(1-x)^{\alpha}(1+x)^{\beta}$, $w=x^{\alpha}e^{-x}$, $w=e^{-x^2}$, the zeros of $p_n^{(\alpha-1,\beta-1)}$, $L_n^{(\alpha-1)}$, $H_n$ respectively (cf. \cite{aho}).

\medskip

As it was pointed out in \cite{ho}, the weighted Fej\'er constants play fundamental role in the energy problem and interpolation as well.

\begin{defi}

The weighted Fej\'er constants on $I$ with respect to $U_n$ and $w$ are:
\begin{equation}C_k:=C_{k,U_n,w}=\frac{\omega^{''}}{\omega^{'}}(u_k)+\frac{w^{'}}{w}(u_k).\end{equation}\end{defi}

\medskip

 We will investigate the local extrema of the energy function (cf. \cite{i}, \cite{ho}). Our main tool will be the differential equation of orthogonal polynomials, and it's transformed version to a Schr\"odinger equation (cf. \cite{j} and \cite{m}). The next general lemma will be useful for the examinations below.

\medskip

\begin{lemma} Let $p_n(x)=\gamma_n\prod_{i=1}^n(x-\zeta_i)$, be a polynomial of degree $n$ with zeros $Z_n=\{\zeta_1, \dots , \zeta_n\}$, for which 
\begin{equation} p_n^{''}(x)+ M_n(x)p_n^{'}(x)+N_n(x)p_n(x)=0.\end{equation}
Let us assume that $M_n$ is the logarithmic derivative of a function $w_n(x)$ which is smooth enough, that is  
$$\left(\log w_n(x)\right)^{'}=M_n(x).$$
Let $T_{w_n}(u_1, \dots , u_{n})$ be the energy function with respect to $w_n$. Then
\begin{equation}\frac{\partial \log  T_{w_n}(u_1, \dots , u_{n})}{\partial u_i}(\zeta_1,\dots ,\zeta_n)=C_{i,w_n,Z_{n}}=0,\end{equation}
\begin{equation}\frac{\partial^2  \log T_{w_n}(u_1, \dots , u_{n})}{\partial u_i\partial u_j}(\zeta_1,\dots ,\zeta_n)=\frac{2}{(\zeta_i-\zeta_j)^2},\end{equation}
and
\begin{equation}\frac{\partial^2 \log  T_{w_n}(u_1, \dots , u_{n})}{\partial u_i^2}(\zeta_1,\dots ,\zeta_n)=-\frac{2}{3}\Phi(\zeta_i),\end{equation}
where $\Phi(x):=\Phi_{w_n}(x)$ is the coefficient of the transformed differential equation:
$$z_n^{''}(x)+\Phi(x)z_n(x)=0,$$
which is satisfied by 
$$z_n(x)=p_n(x)\sqrt{w_n(x)}.$$\end{lemma}

\medskip

\remark
Equality (2) is valid on the domain of definition of the coefficients. $w_n(x)$ is smooth enough, usually means that it is twice differentiable on that domain.

\medskip

\proof 
$$\frac{\partial \log T(u_1, \dots , u_{n})}{\partial u_i}(\zeta_1,\dots ,\zeta_n)=\frac{w_n^{'}}{w_n}(\zeta_i)+2\sum_{1\leq j \leq n\atop j\neq i}\frac{1}{\zeta_i-\zeta_j}$$ $$= \frac{w_n^{'}}{w_n}(\zeta_i)+\frac{p_n^{''}}{p_n^{'}}(\zeta_i)=C_{i,w_n,Z_{n}}=0,$$
where the last equality is ensured by the differential equation. (4) is obvious. To prove (5) at first we have to note that by standard arguments we have  
\begin{equation}\Phi(x) = N_n(x)-\frac{1}{4}M_n^2(x)-\frac{1}{2}M_n^{'}(x),\end{equation}
cf. e.g. \cite{m} Th. 3.6. (The same appears in \cite{q} (5), with $g(x)=x$.)
$$\frac{\partial^2 \log T_{w_n}(u_1, \dots , u_{n})}{\partial u_i^2}=\left(\log w_n\right)^{''}(u_i)-2\sum_{1\leq j \leq n\atop j\neq i}\frac{1}{(u_i-u_j)^2}.$$
Since denoting by $\omega(x)=\omega_{U_n}(x)$
$$\sum_{1\leq j \leq n\atop j\neq i}\frac{1}{(u_i-u_j)^2}=\left(\frac{1}{4}\left(\frac{\omega^{''}}{\omega^{'}}\right)^2-\frac{1}{3}\frac{\omega^{'''}}{\omega^{'}}\right)(u_i),$$
$$\frac{\partial^2 \log  T(u_1, \dots , u_{n})}{\partial u_i^2}(\zeta_1,\dots ,\zeta_n)=\left(\log w_n\right)^{''}(u_i)+\left(-\frac{1}{2}\left(\frac{p_n^{''}}{p_n^{'}}\right)^2+\frac{2}{3}\frac{p_n^{'''}}{p_n^{'}}\right)(\zeta_i).$$
By differentiation of (2) we have
$$\frac{p_n^{'''}}{p_n^{'}}(\zeta_i)=(-M_n^{'}-M_n\frac{p_n^{''}}{p_n^{'}}-N_n)(\zeta_i).$$
Recalling that
$$\frac{p_n^{''}}{p_n^{'}}(\zeta_i)=-M_n(\zeta_i),$$
we get
$$\frac{\partial^2 \log T}{\partial u_i^2}(\zeta_1,\dots ,\zeta_n)=\left(M_n^{'}-\frac{1}{2}M_n^2-\frac{2}{3}M_n^{'}+\frac{2}{3}M_n^2-\frac{2}{3}N_n\right)(\zeta_i)=-\frac{2}{3}\Phi(\zeta_i).$$

\medskip

In the following sections we will investigate the energy function at the zeros of $X_m$-exceptional polynomials. At first we will deal with $X_m$-Laguerre-(I) polynomials, because all the zeros of these polynomials are simple and real. Some of these results can be extended to other families of exceptional polynomials. These results are collected in the last section.

\medskip

\section{$X_m$-Laguerre-(I) Polynomials}

Let $w^{(\alpha)}(x)=x^{\alpha}e^{-x}$, $w^{(\alpha+1)}(x)$ the Laguerre weights on $(0,\infty)$. $L_k^{(\alpha)}(x)$ is the $k^{th}$ Laguerre polynomial of parameter $\alpha$. These classical Laguerre polynomials satisfy the differential equation:
\begin{equation}xy^{''}+(\alpha+1-x)y^{'}+ny=0,\end{equation}
cf. e.g. \cite{sz}, p. 90.
The zeros of $L_k^{(\alpha)}(x)$ are $0<\zeta_{k,1}^{\alpha}< \dots < \zeta_{k,k}^{\alpha}$.

The exceptional Laguerre polynomials of codimension $m\geq 1$ and of the first kind are $\{L_{m,m+n}^{I, (\alpha)}\}_{n=0}^{\infty}$. They are the orthogonal polynomials on $(0, \infty)$ with respect to the weight $\hat{w}^{(\alpha)}_m:=\frac{|x|^{\alpha}e^{-x}}{S^2(x)}$, where $S(x):=S^{(\alpha-1)}_m(x):=L_m^{(\alpha-1)}(-x)$. Furthermore we need the zeros of $L_m^{(\alpha-1)}$, which are $y_1, \dots , y_m$ for simplicity. $L_{m,m+n}^{I, (\alpha)}$ satisfies the differential equation (cf. \cite{gukm}, (75))
\begin{equation}y^{''}(x)+\left(\frac{\alpha+1-x}{x}-\frac{2S^{'}(x)}{S(x)}\right)y^{'}(x)+\left(\frac{m+n}{x}-\frac{\alpha}{x}\frac{2S{'}(x)}{S(x)}\right)y(x)=0.\end{equation}
That is
\begin{equation}M_{m,n}(x)=M_m(x)=\frac{\alpha+1-x}{x}-2\frac{S^{'}(x)}{S(x)}, \ws\ws\ws N_{m,n}(x)=\frac{m+n}{x}-2\frac{S{'}(x)}{S(x)}\frac{\alpha}{x},\end{equation}
cf. (2) and \cite{m} (3.13),(3.14). (8) fulfils on $\mathbb{R}\setminus\{0, -y_1,\dots ,-y_m\}$.

\medskip

\begin{lemma}(\cite{gumm} Prop. 3.2, Cor. 3.1, Prop. 3.4. ) For $\alpha>0$ $L_{m,m+n}^{I, (\alpha)}$ has $m$ simple exceptional zeros in $(-\infty,0)$: $z_{m,n,1}> \dots > z_{m,n,m}$, and $n$ simple regular zeros in $(0, \infty)$: $x_{m,n,1}< \dots < x_{m,n,n}$. The location of these zeros is the following: $0< x_{m,n,1}<\zeta_{n,1}^{\alpha}$, $\zeta_{n-1,j-1}^{\alpha}<x_{m,n,j}<\zeta_{n,j}^{\alpha}$, $ -\zeta_{m,1}^{\alpha}<z_{m,n,1}<0$ and $-\zeta_{m,j}^{\alpha}<z_{m,n,j}< -\zeta_{m-1,j}^{\alpha}$.
Futhermore 
\begin{equation}\lim_{n\to \infty}nx_{m,n,j}=\frac{\left(j_j^{(\alpha)}\right)^2}{4},\end{equation}
where $\{j_j^{(\alpha)}\}_{j\geq 1}$ is the increasing sequence of the positive zeros of the Bessel function $J_{\alpha}(z)$, and the exceptional zeros of $L_{m,m+n}^{I, (\alpha)}$ converge to the $m$ zeros of $L_m^{(\alpha-1)}(-x)$. \end{lemma}

\medskip 

For simplicity let us denote by $z_i:=z_{m,n,i}$ and by $x_j:= x_{m,n,j}$. That is 
\begin{equation} L_{m.m+n}^{I, (\alpha)}= \left(\frac{1}{m!}\prod_{i=1}^m(x-z_i)\right) \frac{(-1)^n}{n!}\prod_{j=1}^n(x-x_j)=P_{m,n}(x)q_{m,n}(x).\end{equation}

\medskip

\remark
By (21) of \cite{gumm} the leading coefficient of $L_{m,m+n}^{I, (\alpha)}$ is $\frac{(-1)^n}{m!n!}$. Also by (21) of \cite{gumm} $L_{m,m+n}^{I, (\alpha)}(0)\neq 0$. 

\medskip

 With the notations of Lemma 1 and by (8) $w_n=\hat{w}^{(\alpha+1)}_m$, which doesn't depend on $n$. So the energy function is
 $$T_{\hat{w}^{(\alpha+1)}_m}\left(u_1, \dots ,u_{m+n}\right)=\prod_{i=1}^{m+n}\hat{w}^{(\alpha+1)}_m(u_i)\prod_{1\leq i < j \leq n+m}(u_i-u_j)^2.$$
 The following statement is the extension of Th. 2 of \cite{dy} to $X_m$-Laguerre(I) polynomials.

\medskip

\begin{theorem} Let $\alpha \geq 1$ and $n \geq 0$. Then the logarithmic energy function, $ \log T_{\hat{w}^{(\alpha+1)}_m}(u_1, \dots , u_{m+n})$ has a saddle point at \\ $Z_{m,n}=\{z_{m,n,1}, \dots , z_{m,n,m}, x_{m,n,1}, \dots , x_{m,n,n}\}$, which are the zeros of $L_{m,m+n}^{I, (\alpha)}$. More precisely
$\frac{\partial^2 \log  T_{\hat{w}^{(\alpha+1)}_m}(u_1, \dots , u_{n})}{\partial u_i^2}\left(Z_{m,n}\right)$ is positive if $u_i$ is one of the first $m$ variables, and it is negative if  $u_i$ is one of the last $n$ variables.
\end{theorem}

\proof According to (3), the first partial derivatives of the logarithmic energy function are zero at the zeros of $L_{m,m+n}^{I, (\alpha)}$.
 Also by Lemma 1, the Hesse matrix of $\log T_{\hat{w}^{(\alpha+1)}_m}$ (denoted by $H$) is the folowing:
$$H_{i,j}=\frac{2}{(u_i-u_j)^2}, \ws \ws \mbox{and } \ws\ws H_{i,i}=-\frac{2}{3}\Phi(u_i).$$
If $u_i=x_i$ a positive zero of $L_{m,m+n}^{I, (\alpha)}$, then
$$H_{i,i}=-\frac{\alpha+1}{x_i^2}+2\sum_{l=1}^m\left(\frac{1}{(x_i+y_l)^2}-\frac{1}{(x_i+\xi_l)^2}\right)-2\sum_{1\leq j \leq n \atop j\neq i}\frac{1}{(x_i-x_j)^2}.$$
where $\xi_j=|z_j|$. For $j=1, \dots , m$ we have
$$0<\xi_j<\zeta_{m,j}^{\alpha}<\zeta_{m+1,j+1}^{\alpha-1}<\zeta_{m,j+1}^{\alpha-1}=y_{j+1},$$
where  the second inequality is ensured by Lemma 2, the third one fulfils because $L_m^{(\alpha)}=-\left(L_{m+1}^{(\alpha-1)}\right)^{'}$, and the last one by the interlacing property of the zeros of orthogonal polynomials. Rearranging the sum above we have
$$H_{i,i}=-\frac{\alpha+1}{x_i^2}+ \frac{2}{(x_i+y_1)^2}+2\sum_{l=1}^{m-1}\left(-\frac{1}{(x_i+\xi_l)^2}+\frac{1}{(x_i+y_{l+1})^2}\right)$$ $$-\frac{2}{(x_i+\xi_m)^2}-2\sum_{1\leq j \leq n \atop j\neq i}\frac{1}{(x_i-x_j)^2},$$
which is negative if $\alpha\geq 1$ for all $n\geq 0$.

If $u_i=z_i$ is a negative zero, we will show that \\ $-\Phi(z_i)=\left(\frac{1}{2}M_n^{'}+\frac{1}{4}M_n^2-N_n\right)(z_i)>0$.
Since
$$M_n^{'}(x)=-\frac{\alpha+1}{x^2}-2\left(\frac{S^{''}}{S}-\left(\frac{S^{'}}{S}\right)^2\right)(x),$$
we have to compute $\frac{S^{''}}{S}$. Recalling that $S(x)=L_m^{(\alpha-1)}(-x)$, by (7)
$$\frac{S^{''}}{S}(x)=-\frac{\alpha+x}{x}\frac{S^{'}}{S}(x)+\frac{m}{x}.$$ 
Substituting this we have
$$-\Phi(x)=\left(-\frac{\alpha+1}{2x^2}+\frac{\alpha+x}{x}\frac{S^{'}}{S}(x)-\frac{m}{x}+\left(\frac{S^{'}}{S}\right)^2(x)\right)$$ $$+\frac{1}{4}\left(\frac{\alpha+1}{x}-1-2\frac{S^{'}}{S}(x)\right)^2-\left(\frac{m+n}{x}-2\frac{\alpha}{x}\frac{S^{'}}{S}(x)\right)$$
$$=2\left(\frac{S^{'}}{S}(x)+\frac{1}{2}+\frac{2\alpha-1}{4x}\right)^2-\frac{1}{4}-\frac{2\alpha^2-4\alpha+3}{8x^2}-\frac{4m+2n+3\alpha}{2x}.$$
So
$$-\Phi(-\xi_i)=2(\cdot)^2+f(\xi_i),$$
where
$$f(\xi_i)=\frac{-2\xi_i^2+ 4(4m+2n+3\alpha)\xi_i -(2\alpha^2-4\alpha+3)}{8\xi_i^2}.$$
That is $f(\xi_i)$ is positive if $\xi_i$ is between the two zeros of the numerator: $m_1\leq \xi_i\leq m_2$.
By Lemma 2 $\xi_m\leq\zeta_{m,m}^{\alpha}\leq 2m+\alpha+1+\sqrt{(2m+\alpha+1)^2+\frac{1}{4}-\alpha^2}$ (cf. \cite{sz}, Th. 6.31.2.). It can be seen that if $\alpha\geq 1$, for all $n\geq 0$
$$2m+\alpha+1+\sqrt{(2m+\alpha+1)^2+\frac{1}{4}-\alpha^2}$$ $$\leq 4m+2n+3\alpha+\sqrt{(4m+2n+3\alpha)^2-(\alpha-1)^2-\frac{1}{2}}=m_2.$$
Also by Lemma 2 and by \cite{sz} Th. 6.31.3 
$$\xi_1>\zeta_{m,1}^{(\alpha)}>\frac{\left(\frac{j_1^{(\alpha)}}{2}\right)^2}{m+\frac{\alpha+1}{2}}>\frac{(2\alpha+2\pi-1)^2}{16m+8(\alpha+1)},$$
where $j_1^{(\alpha)}$ is the first positive zero of the Bessel function with parameter $\alpha$, and see e. g. \cite{e} for the last estimation.
As above, it can be seen that for all $n\geq 0$
$$m_1=\frac{(\alpha-1)^2+\frac{1}{2}}{4m+2n+3\alpha+\sqrt{(4m+2n+3\alpha)^2-(\alpha-1)^2-\frac{1}{2}}}\leq\frac{\left(\alpha+\pi-\frac{1}{2}\right)^2}{4m+2(\alpha+1)},$$
that is
$$-\Phi(-\xi_i)>0, \ws \ws \ws i=1, \dots ,m \ws \ws \ws \forall n\geq 0.$$
Since the Hessian is diagonally dominant, the result is proved.

\medskip

Summarizing in the previous theorem we had $w_n(x)=\frac{|x|^{\alpha+1}}{\left(L_m^{(\alpha-1)}(-x)\right)^2}$ which depends on $m$ but independent of $n$, and it is positive on  the interval of orthogonality $I=[0,\infty)$, but $S$ has $m$ zeros away of $I$. The solution of the differential equation is $P_{m,n}q_{m,n}$, where $P_{m,n}$ is a polynomial of degree $m$ with zeros $\{z_1, \dots ,z_m\}$ away of $I$, and $q_{m,n}$ is a polynomial of degree $n$ with zeros $\{x_1, \dots ,x_n\}$ in $I$. The logarithmic energy function with respect to $w_n$ with $n+m$ variables $\{u_1, \dots ,u_{n+m}\}$ has a saddle point at the zeros of $L_{m.m+n}^{I, (\alpha)}$. We are looking for an energy function wich has a maximum at some zeros of $L_{m.m+n}^{I, (\alpha)}$. To this purpose let us denote by 

$$v:=v^{(\alpha+1)}_{m,n}:=\hat{w}^{(\alpha+1)}_mP_{m,n}^2
=\frac{w_1P^2}{S^2}.$$ 
That is $v$ is a new weight function on $(0,\infty)$, which depends on $n$.

and $T_v(u_1, \dots , u_n)=\prod_{j=1}^nv(u_j)\prod_{1\leq i<j\leq n}(u_i-u_j)^2.$ With these notations
$$\frac{\partial \log  T_{\hat{w}^{(\alpha+1)}_m}(u_1, \dots , u_{n+m})}{\partial u_{m+i}}(z_1,\dots ,x_n)=\frac{\left(\hat{w}^{(\alpha+1)}_m\right)^{'}}{\hat{w}^{(\alpha+1)}_m}(x_i)+2\frac{P^{'}}{P}(x_i)$$ $$+2\sum_{1\leq j\leq n \atop j\neq i}\frac{1}{x_i-x_j}=\frac{v^{'}}{v}(x_i)+2\sum_{1\leq j\leq n \atop j\neq i}\frac{1}{x_i-x_j}=\frac{\partial \log  T_v(u_1, \dots , u_{n})}{\partial u_{i}}(x_1,\dots ,x_n),$$
and
$$\frac{\partial^2 \log  T_{\hat{w}^{(\alpha+1)}_m}(u_1, \dots , u_{n+m})}{\partial u_{m+i}^2}(z_1,\dots ,x_n)=\left(\log v\right)^{''}(x_i)-2\sum_{1\leq j\leq n \atop j\neq i}\frac{1}{(x_i-x_j)^2}$$ $$=\frac{\partial^2 \log  T_v(u_1, \dots , u_{n})}{\partial u_{i}^2}(x_1,\dots ,x_n).$$
The computations above show that it is not a suprising idea taking into consideration only the regular zeros of the exceptional polynomials, and using the new weight. Now we have the following differential equation for $q_{m,n}=q$
$$q^{''}(x)+M_{1,n}(x)q^{'}(x)+N_{1,n}(x)q(x)=0,$$
where all the expressions depend on $m$ too, and
$$M_{1,n}(x)=M_n(x)+2\frac{P_{m,n}^{'}}{P_{m,n}}(x), \ws\ws N_{1,n}(x)=N_n(x)+\frac{P_{m,n}^{''}}{P_{m,n}}(x)+M_n(x)\frac{P_{m,n}^{'}}{P_{m,n}}(x).$$
Obviously $M_{1,n}(x)=\left(\log v\right)^{'}(x)$.  Moreover $q_{m,n}(x)\sqrt{v(x)}$ fulfils the differential equation:
$$f^{''}(x)+\Phi_1(x)f(x)=0, \ws \ws \ws \mbox{where }\ws\ws \Phi_1(x)=N_{1,n}(x)-\frac{1}{4}M_{1,n}^2(x)-\frac{1}{2}M_{1,n}^{'}(x). $$
So as a corollary of the previous theorem we can state

\medskip

\begin{theorem} Let $\alpha \geq 1$. The positive zeros of $L_{m,m+n}^{I, (\alpha)}$: \\$X_{m,n}=\{x_{m,n,1}, \dots , x_{m,n,n}\}$ is the unique set of minimal energy (or Fekete set) with respect to the external field represented by the weight $\left(v^{(\alpha+1)}_{m,n}\right)^{\frac{1}{2(n-1)}}$.\end{theorem}

\proof Denoting by $v=v^{(\alpha+1)}_{m,n}$ and by $P=P_{m,n}$, $q=q_{m,n}$, and the energy function is $T_v$ on $[0,\infty)^n$.
According to the computations above
$$\frac{\partial \log T_v(u_1, \dots , u_n)}{\partial u_i}(x_1, \dots , x_n)=C_{i,v,X}$$ $$=\frac{\left(\hat{w}^{(\alpha+1)}_m\right)^{'}}{\hat{w}^{(\alpha+1)}_m}(x_i)+\frac{\left(L_{m,m+n}^{I, (\alpha)}\right)^{''}}{\left(L_{m,m+n}^{I, (\alpha)}\right)^{'}}(x_i)=C_{m+i,\hat{w}^{(\alpha+1)}_m,Z}=0.$$
Similarly, the computations above entails that at $(x_1,\dots ,x_n)$ the energy funtion has a local maximum. Computing the Hesse matrix $H$ of $\log T$ at any point of $(0,\infty)^n$, we can investigate the global behavior of the energy function. By Lemma 1 we have
$$H_{i,j}=\frac{2}{(u_i-u_j)^2},$$
and
$$H_{i,i}=\left(\log v\right)^{''}(u_i)-2\sum_{1\leq j \leq n\atop j\neq i}\frac{1}{(u_i-u_j)^2}.$$
$$\left(\log v\right)^{''}(u_i)=2\left(\frac{P^{'}}{P}-\frac{S^{'}}{S}\right)^{'}(u_i)-\frac{\alpha+1}{u_i^2}$$
$$= 2\sum_{j=1}^m\left(\frac{1}{(u_i+y_j)^2}-\frac{1}{(u_i+\xi_j)^2}\right)-\frac{\alpha+1}{u_i^2}.$$
Proceeding as above  since $\alpha \geq 1$, it can be seen that for any $U_n \subset (0,\infty)$ $H_{i,i}<0$ for $i=1, \dots, n$, for all $n\geq 0$. So $-H$ is real, symmetric, strictly diagonally dominant and positive definite for all $U_n\subset (0,\infty)$. Furthermore since $\alpha>0$, $T$ tends to zero at the boundary of the domain, so it has a unique maximum at $X$.

\medskip

\remark 

\noindent (1) As it was pointed out in \cite{ho}, the assumption: $-\left(\log v\right)^{''}>0$ on $(0,\infty)$ ensures the uniqueness of the system of minimal energy. In Lemma 1 we have seen that this is exactly the opposite of the first term in $H_{i,i}$ that is 
$$-\left(\log v\right)^{''}=-2\sum_{j=1}^{m-1}\left(-\frac{1}{(u_i+\xi_j)^2}+\frac{1}{(u_i+y_{j+1})^2}\right)$$ 
$$- \frac{2}{(u_i+y_1)^2}+\frac{2}{(u_i+\xi_m)^2}+\frac{\alpha+1}{u_i^2},$$
which is positive by the previous calculation. 

\noindent (2) If $0<\alpha<1$ and $i\neq 1$, then $\frac{2}{(u_i+y_1)^2}-\frac{2}{(u_i-u_1)^2}<0$, so for $i > 1$,  $H_{i.i}<0$ for all $U_n \subset (0,\infty)$, but we can choose a $U_n$ such that $H_{1,1}>0$.

\medskip

Now we turn to the Egerv\'ary-Tur\'an interpolatory problem, which is to find an interpolatory process of the lowest degree and of the smallest norm. A first solution can be found in \cite{jo}, a result to Markov-Sonin polynomials in \cite{b} and the generalization to any classical weights in \cite{rs}, to general weights in \cite{ho}.
Below we denote by $\hat{l}_k(x)$ any polynomials of arbitrary degree for which $\hat{l}_k(x_i)=\delta_{ki}, \ws i=1, \dots , n$.
\begin{defi} Let $w$ be a positive weight. The interpolatory system of polynomials $\hat{l}_k(x), \ws k=1, \dots , n$ is $w$-stable on $(a,b)$ if for all $y_1, \dots , y_n \geq 0$
$$ 0\leq w(x)\sum_{k=1}^n\frac{\hat{l}_k(x)}{w(x_k)}y_k \leq \max_ky_k, \ws\ws x\in (a,b).$$
A $w$-stable interpolatory system on $(a,b)$ is most economical, if 
$$ \sum_{k=1}^n \deg \left(\hat{l}_k(x)\right)$$
is minimal.\end{defi}

\medskip

With the notations of Lemma 1, the Hermite interpolatory operator on the zeros of $p_n$ of the function $\frac{1}{w_n}$ coincides with the Gr\"unwald operator of it:
with the usual notation $l_k(x)=\frac{p_n(x)}{p_n^{'}(\zeta_k)(x-\zeta_k)}$, by the differential equation (2),
$$H_n\left(\frac{1}{w_n},\left(\frac{1}{w_n}\right)^{'},x\right)=\sum_{k=1}^n\left(1-\frac{p_n^{''}}{p_n^{'}}(\zeta_k)(x-\zeta_k)\right)l_k^2(x)\frac{1}{w_n(\zeta_k)}$$
$$+\sum_{k=1}^n(x-\zeta_k)l_k^2(x)\left(\frac{1}{w_n}\right)^{'}(\zeta_k)=\sum_{k=1}^nl_k^2(x)\frac{1}{w_n(\zeta_k)}=:Y_n\left(\frac{1}{w_n},x\right),$$
and on the righthand-side there is exactly the  Gr\"unwald operator with respect to  $\frac{1}{w_n}$. So if  $\left(\frac{1}{w_n}\right)^{(2n)}\geq 0$ for all $n\geq 0$, then the error formula of Hermite interpolation ensures that the Gr\"unwald operator is $w_n$-stable and it is the most economical as well. Unfortunately the assumption on the even derivatives does not fulfil when $w_n=\hat{w}^{(\alpha+1)}_m$ and $x\in (z_m, x_n)$. But we can state the following

\medskip

\begin{theorem} If $\alpha>1$, the Gr\"unwald operator on the positive zeros of $L_{m,m+n}^{I, (\alpha)}$ is $v^{(\alpha+1)}_{m,n}$-stable and it is also the most economical.\end{theorem}

\medskip

\begin{lemma} If $\alpha>1$, 

$$\left(\frac{1}{v^{(\alpha+1)}_{m,n}}\right)^{(2n)}(x)>0, \ws\ws\ws x \in (0,\infty), \ws\ws\ws n=0,1,2,\dots .$$\end{lemma}

\proof  With the notations above,
$$\frac{1}{v^{(\alpha+1)}_{m,n}}=\frac{S^2(x)e^x}{P^2(x)x^{\alpha+1}}=\frac{e^x}{x^{\alpha-1}}\left(1+\frac{y_1}{x}\right)^2\frac{1}{(x+\xi_m)^2}\prod_{i=1}^{m-1}\left(1+\frac{y_{i+1}-\xi_i}{x+\xi_i}\right)^2$$
$$=\left(1+2\frac{y_1}{x}+\frac{y_1^2}{x^2}\right)\frac{e^x}{x^{\alpha-1}}\sum_{k_1=0,1,2, \dots \atop k_{m-1}=0,1,2}\frac{c_{k_1, \dots ,k_{m-1}}}{(x+\xi_1)^{k_1}\dots(x+\xi_{m-1})^{k_{m-1}}(x+\xi_m)^2},$$
where all the coefficients $c_{k_1, \dots ,k_{m-1}}$ are positive. Expressing these products as products of Laplace transforms, let
$$g(\alpha-1)=\sum_{k_1=0,1,2, \dots \atop k_{m-1}=0,1,2}\frac{c_{k_1, \dots ,k_{m-1}}}{\Gamma(\alpha-1)(k_1-1)^{*}!\dots (k_{m-1}-1)^{*}!}$$ 
$$ \times\int_0^{\infty}\dots \int_0^{\infty}t_0^{\alpha-2}t_1^{k_1-1}\dots t_{m-1}^{k_{m-1}-1}t_me^{-(t_1\xi_1+ \dots +t_m\xi_m)^{*}}e^{-x(t_0-1+t_1+\dots +t_m)^{*}}dt_0^{*} \dots dt_{m}^{*},$$ 
where $(\cdot)^{*}$ means that if $k_i=0$ then the $k_i^{th}$ term is missing, and the integral against the $k_i^{th}$ variable is also missing.
So
$$\frac{1}{v^{(\alpha+1)}_{m,n}}=g(\alpha-1)+2y_1g(\alpha)+y_1^2g(\alpha+1),$$
and this implies that the even derivatives of the reciprocal of the weight function are all positive.

\medskip

\proof ( of Th. 3)
First of all, let us observe that for $k=1, \dots ,n$
$$l_{m,k}(x):=\frac{(Pq)(x)}{(Pq)^{'}(x_k)(x-x_k)}=\frac{P(x)}{P(x_k)}\frac{q(x)}{q^{'}(x_k)(x-x_k)}=:\frac{P(x)}{P(x_k)}l_k(x),$$
where we used the notations of Th. 2. Thus as above
$$\hat{w}^{(\alpha+1)}_m(x)\left(\sum_{k=1}^n\left(1-\frac{\left(L_{m,m+n}^{I,(\alpha)}\right)^{''}}{\left(L_{m,m+n}^{I,(\alpha)}\right)^{'}}(x_k)(x-x_k)\right)l_{m,k}^2(x)\frac{1}{\hat{w}^{(\alpha+1)}_m(x_k)}\right)$$ $$+\hat{w}^{(\alpha+1)}_m(x)\left(\sum_{k=1}^n(x-x_k)l_{m,k}^2(x)\left(\frac{1}{\hat{w}^{(\alpha+1)}_m}\right)^{'}(x_k)\right)$$ $$=\hat{w}^{(\alpha+1)}_mP^2(x)\sum_{k=1}^nl_k^2(x)\frac{1}{\left(\hat{w}^{(\alpha+1)}_mP^2\right)(x_k)}=v(x)\sum_{k=1}^n\frac{l_k^2(x)}{v(x_k)}.$$
Hence $\sum_{k=1}^n\frac{l_k^2(x)}{v(x_k)}$ is the Hermite interpolatory polynomial of $\frac{1}{v}$ on $x_1, \dots ,x_n$, and by the error formula and Lemma 3
$$1-v(x)\sum_{k=1}^n\frac{l_k^2(x)}{v(x_k)}>0 \ws\ws\ws x>0,$$
which ensures that norm of the Gr\"unwald operator is at most one, so it is $v$-stable, and by the construction it is the most economical.

\medskip
      
Recalling the notion of the energy function $T=T_v$ on $(0,\infty)^n$ with respect to $v=v^{(\alpha+1)}_{m,n}$, we can write 
$$-\log\left( (T_v)^{\frac{1}{n(n-1)}}\right)=\frac{2}{n(n-1)}\sum_{1\leq i<j\leq n}k_n(u_i,u_j)=:d(u_1,\dots ,u_n),$$
where
$$k_n(x,y):=-\log \left(\frac{c}{n}|x-y|v^{\frac{1}{2(n-1)}}(x)v^{\frac{1}{2(n-1)}}(y)\right),$$
is the modified kernel (with $\frac{c}{n}$), which is symmetric, lower semicontinuous and positive. Here, taking into consideration that by Lemma 2 and by (11) $\frac{P^2(x)}{S^2(x)}$ can be estimated on $(0,\infty)$ by a constant, $c$ is an absolute constant.

If the kernel is independent of $n$ it is proved that the $n^{th}$ transfinite diameter,
$$\delta_n:=\inf_{u_1, \dots ,u_n \in H}\frac{2}{n(n-1)}\sum_{1\leq i<j\leq n}k(u_i,u_j),$$ 
has a limit with $n \to \infty$, cf. \cite{fana}. This limit is called as the transfinite diameter of $H$ with respect to the kernel. It is also obvious, that the lower semicontinuity of the kernel ensures that on compact sets the infinum is a minimum. Although the positive axis is not compact, but the weight tends to zero quickly enough at the endpoints, so here we also have a minimum, as it was pointed out in Th. 2. Next we will prove that the $n^{th}$ transfinite diameter, has a limit in our case as well. When the kernel is independent of $n$,  usually $\delta_n$ is an increasing sequence.  Let us denote by
$$d_n=\inf_{u_1, \dots ,u_n \in (0,\infty)}d(u_1,\dots ,u_n).$$
In our case, by these standard arguments it can be shown only that $d_n\geq d_{n-1}-O\left(\frac{1}{n}\right)$, which is not enough. In Th. 2. we saw that $d_n=d(x_1, \dots ,x_n)$, where the extremal set is the set of the zeros of $L_{m,m+n}^{I, (\alpha)}$. In order to prove that $d_n$ has a limit we will estimate the difference of two consecutive terms of the sequence.

We have to mention here that the existence of the limit can be proved by potential theoretic tools as well. As it was appeard at first in \cite{ms}, defining $\sigma_n(x):=v(a_nx)^{\frac{1}{2n}}$, where $a_n:=a_{m,n}$ is the Mhaskar-Rahmanov-Saff number with respect to $\sqrt{v}$, where $v=v^{(\alpha+1)}_{m,n}$, $\{\sigma_n(x)\}_n$ is a sequence of weights on the common compact support $[0,1]$. Since $\|\sqrt{v(x)}p_n(x)\|_{\infty,[0,\infty)}=\|\sqrt{v(x)}p_n(x)\|_{\infty,[0,a_n)}$, one can find that if $\sigma_n$ converges well on $[0,1]$ to a weight $\sigma$, then the Chebyshev constant defined by $t_n(\sigma_n)^{\frac{1}{n}}:=\left(\inf_{p_n(x)=x^n+\dots}\|\sigma_n(x)^np_n(x)\|_{\infty,[0,1]}\right)^{\frac{1}{n}}$ has similar asymptotic behavior as $t_n(\sigma)^{\frac{1}{n}}$ has. If $\sigma$ is nice (admissible) this limit exists and it coincides with $\lim_{n\to \infty}\frac{\left(\inf_{p_n(x)=x^n+\dots}\|\sqrt{v(x)}p_n(x)\|_{\infty,[0,\infty)}\right)^{\frac{1}{n}}}{a_n}$. Let us obseve that $a_n\left(x^{\alpha+1}e^{-x}\right)\sim n$ and by Lemma 2, $a_n(\sqrt{v})\sim n$ (cf. the normalization above). So by the replacement $x=a_n\xi$, $y=a_n\eta$, $k_{n,v}(x,y)=-\log \left(\frac{|x-y|}{a_n}v^{\frac{1}{2(n-1)}}(x)v^{\frac{1}{2(n-1)}}(y)\right)=-\log \left(|\xi-\eta|v^{\frac{1}{2(n-1)}}(a_n\xi)v^{\frac{1}{2(n-1)}}(a_n\eta)\right)=-\log \left(|\xi-\eta|w_n(\xi)w_n(\eta)\right)+ e(\xi,\eta)$ is a kernel function on $[0,1]$ with varying weights (cf. \cite{t}), where $e(\xi,\eta)$ is an error term. Thus the corresponding expression is (cf.\cite{fana})\\ $M_n:=\sup_{U_n\subset(0,a_n)}\inf_{x\in (0,a_n)}\frac{1}{n}\sum_{i=1}^nk_n(x,u_i)$\\$=-\log \left(\inf_{U_n\subset(0,a_n)}\sup_{x\in (0,a_n)}\frac{\left(\sqrt{v(x)}p_n(x)\right)^{\frac{1}{n}}}{a_n}v(x)^{\frac{1}{2n(n-1)}}\left(\prod_{i=1}^nv(u_i)\right)^{\frac{1}{n}}\right)$, where the $n^{th}$ Chebyshev constant appears, then a term which has to tend to zero, and a correction term from the integral of $-\log v$. (cf. \cite{st} p. 162.) Finally taking into consideration the connection between the Chebyshev constant and the transfinite diameter, one can get the fact of convergence by this way too. Here we give an estimation on the speed of convergence by classical methods.

\medskip

\begin{theorem} Let $\alpha \geq 1$. With the notations above $d_n=d_{n-1}+O\left(\frac{\log^2n}{n^2}\right)$.
\end{theorem}

\proof  For simplicity, as above $x_1,\dots ,x_n$, $\nu_1,\dots ,\nu_{n-1}$ are the positive zeros of $L_{m,m+n}^{I, (\alpha)}$ and $L_{m,m+n-1}^{I, (\alpha)}$; and $\zeta_{1,n}, \dots ,\zeta_{n,n}$, $\zeta_{1,n-1}, \dots ,\zeta_{n-1,n-1}$ are the zeros of $L^{\alpha}_n$ and $L^{\alpha}_{n-1}$ respectively. Now
$$d_n=\frac{2}{n(n-1)}\sum_{1\leq i < j \leq n}k_n(x_i,x_j)=\frac{1}{n}\sum_{l=1}^n\frac{2}{(n-1)(n-2)}\sum_{1\leq i < j \leq n \atop i,j \neq l}k_n(x_i,x_j)$$
$$=\frac{1}{n}\sum_{l=1}^n\frac{2}{(n-1)(n-2)}\sum_{1\leq i < j \leq n \atop i,j \neq l}k_{n-1}(x_i,x_j)+\log\frac{n}{n-1}$$ $$+\frac{1}{n}\sum_{l=1}^n\frac{2}{(n-1)(n-2)^2}\sum_{1\leq i < j \leq n \atop i,j \neq l}\log\left|\frac{P_{n-1}(x_i)P_{n-1}(x_j)}{P_{n}(x_i)P_{n}(x_j)}\right|$$ $$+\frac{1}{n}\sum_{l=1}^n\frac{2}{(n-1)^2(n-2)^2}\sum_{1\leq i < j \leq n \atop i,j \neq l}\log\left|\frac{P_{n}(x_i)P_{n}(x_j)}{S(x_i)S(x_j)}\right|$$ $$+\frac{1}{n}\sum_{l=1}^n\frac{1}{(n-1)^2(n-2)^2}\sum_{1\leq i < j \leq n \atop i,j \neq l}\left((\alpha+1)\log x_i-x_i+(\alpha+1)\log x_j-x_j\right).$$
Let us denote by
$$x_{i,l}=\left\{\begin{array}{ll}x_i, \ws\ws  1\leq i\leq l-1\\
x_{i+1}, \ws\ws  l\leq i \leq n-1\end{array}\right. .$$
Taking into account that
$$k_{n-1}(x_{i,l},x_{j,l})=k_{n-1}(\nu_i,\nu_j)-\log\left|\frac{x_{i,l}-x_{j,l}}{\nu_i-\nu_j}\right|-\frac{1}{n-2}\log\left|\frac{\frac{P_{n-1}}{S}(x_{i,l})\frac{P_{n-1}}{S}(x_{j,l})}{\frac{P_{n-1}}{S}(\nu_{i})\frac{P_{n-1}}{S}(\nu_{j})}\right|$$ $$-\frac{\alpha+1}{2(n-2)}\left(\log\frac{x_{i,l}}{\nu_i}+\log\frac{x_{j,l}}{\nu_j}\right)+\frac{1}{2(n-2)}\left(x_{i,l}+x_{j,l}-\nu_i-\nu_j\right),$$
and denoting by $\log V(U_n):=\sum_{1\leq i < j \leq n}\log|u_j-u_i|$ we have
$$d_n=d_{n-1}$$ 
$$+\left(\log\frac{n}{n-1}+\frac{2}{(n-1)(n-2)}\log V(\Xi_{n-1})-\frac{2}{n(n-1)}\log V(\Xi_{n})\right)$$ 
$$+\left(\frac{2}{(n-1)(n-2)}\left(\log V(X_{n-1})-\log V(\Xi_{n-1})\right)\right.$$ $$\left.-\frac{2}{n(n-1)}\left(\log V(X_{n})-\log V(\Xi_{n})\right)\right)$$
$$\frac{1}{n}\sum_{l=1}^n\frac{2}{(n-1)(n-2)^2}\sum_{1\leq i < j \leq n \atop i,j \neq l}\log\left|\frac{P_{n-1}(x_i)P_{n-1}(x_j)}{P_{n}(x_i)P_{n}(x_j)}\right|$$
$$+\frac{1}{n}\sum_{l=1}^n\frac{2}{(n-1)(n-2)^2}$$ $$\times\sum_{1\leq i < j \leq n-1}\left(\log\left|\frac{\frac{P_{n-1}}{S}(\nu_{i})\frac{P_{n-1}}{S}(\nu_{j})}{\frac{P_{n-1}}{S}(x_{i,l})\frac{P_{n-1}}{S}(x_{j,l})}\right|+\frac{1}{n-1}\log\left|\frac{P_{n}}{S}(x_{i,l})\frac{P_{n}}{S}(x_{j,l})\right|\right)$$
$$+\frac{1}{n}\sum_{l=1}^n\frac{\alpha+1}{(n-1)(n-2)^2}\sum_{1\leq i < j \leq n-1}\left(\log\frac{\nu_i\nu_j}{x_{i,l}x_{j,l}}+\frac{1}{n-1}\log(x_{i,l}x_{j,l})\right)$$
$$+\frac{1}{n}\sum_{l=1}^n\frac{1}{(n-1)(n-2)^2}\sum_{1\leq i < j \leq n-1}\left(x_{i,l} + x_{j,l}-\nu_i- \nu_j-\frac{1}{n-1}(x_{i,l} + x_{j,l})\right)$$
$$=d_{n-1}+M_1+M_2+M_3+M_4+M_5+M_6.$$
We estimate the error terms $M_i$ below.

\noindent $M_3=O\left(\frac{1}{n^{\frac{3}{2}}}\right)$. Indeed, recalling that by (11) and the remark after it, the leading coefficient of $P_n$ is $\frac{1}{m!}$,
$$|M_3|\leq \frac{2}{n(n-2)}\sum_{i=1}^n\left|\log\left|\frac{P_{n-1}}{P_n}(x_i)\right|\right|$$ $$\leq \sum_{j=1}^m\frac{2}{n(n-2)}\sum_{i=1}^n\left|\log\left(1+\frac{z_{m,n,j}-z_{m,n-1,j}}{x_i-z_{m,n,j}}\right)\right|.$$
By Lemma 2, $\Delta:=|x_1-z_{m,n,j}|\geq\frac{c}{m}$, and $|z_{m,n,j}-z_{m,n-1,j}|\leq \varepsilon_n$, and $\Delta_l=x_{l+1}-x_l>\frac{1}{4}(\left(j_{l+1}^{(\alpha)}\right)^2-\left(j_l^{(\alpha)}\right)^2-\frac{2\varepsilon_n}{n}>c\frac{l}{n}$ (cf. \cite{e} (1.6) and 1.4). That is 
$$|M_3|\leq \frac{\varepsilon_n}{n^2}c(m)\sum_{i=1}^n\frac{1}{\Delta+\sum_{l=1}^{i-1}\Delta_l}\leq \frac{c(m)}{n}\sum_{i=1}^n\frac{1}{\frac{n}{m}+i^2}=O\left(\frac{1}{n^{\frac{3}{2}}}\right).$$

\noindent $M_4=O\left(\frac{1}{n^{\frac{3}{2}}}\right)$. As previously
$$\frac{1}{n}\sum_{l=1}^n\frac{2}{(n-1)(n-2)^2}$$ $$\times\sum_{1\leq i < j \leq n-1}\left|\log\left|\frac{P_{n-1}}{S}(\nu_i)\frac{P_{n-1}}{S}(\nu_j)\right|-\log\left|\frac{P_{n-1}}{S}(x_{i,l})\frac{P_{n-1}}{S}(x_{j,l})\right|\right|$$ 
$$+\frac{1}{n}\sum_{l=1}^n\frac{2}{(n-1)^2(n-2)^2}\sum_{1\leq i < j \leq n-1}\left|\log\left|\frac{P_{n-1}}{S}(x_{i,l})\frac{P_{n-1}}{S}(x_{j,l})\right|\right|$$ 
$$ \leq\frac{2}{(n-1)(n-2)}\sum_{j=1}^m\sum_{i=1}^{n-1}\left|\log\frac{\nu_i-z_{m,n-1,j}}{\nu_i+y_j}\right|+\frac{2}{n(n-1)}\sum_{j=1}^m\sum_{i=1}^{n}\left|\log\frac{x_{i}-z_{m,n-1,j}}{x_{i}+y_j}\right|$$ $$
+\frac{c(m)}{n^2}+\frac{2}{n(n-1)^2}\sum_{j=1}^m\sum_{i=1}^{n}\left|\log\frac{x_{i}-z_{m,n-1,j}}{x_{i}+y_j}\right|,$$
and the computation can be finished as above.

\noindent $M_5=O\left(\frac{\log n}{n^{2}}\right)$. 
$$M_5=\frac{1}{(n-1)(n-2)}\sum_{i=1}^{n-1}\log\nu_i-\frac{1}{n(n-1)}\sum_{i=1}^{n}\log x_i=-\frac{\log x_n}{n(n-1)}$$ $$+\frac{1}{(n-1)(n-2)}\sum_{i=1}^{n-1}\log\frac{\nu_i}{x_i}+\frac{2}{n(n-1)(n-2)}\sum_{i=1}^{n}\log x_i $$ $$=O\left(\frac{\log n}{n^{2}}\right)+\frac{1}{(n-1)(n-2)}\sum_{i=1}^{n-1}\log\left(1+\frac{\nu_i-x_i}{x_i}\right)=O\left(\frac{\log n}{n^{2}}\right),$$
where the last estimation is ensured by Lemma 2, namely $\nu_i-x_i\leq c\Delta_i$.

\noindent $M_2=O\left(\frac{\log^2 n}{n^{2}}\right)$. Indeed let us compute for instance
$$\frac{2}{n(n-1)}\sum_{1\leq i < j \leq n}\log\frac{x_j-x_i}{\zeta_{j,n}-\zeta_{i,n}}\leq\frac{2\varepsilon_n}{n^2(n-1)}\sum_{1\leq i < j \leq n}\frac{1}{\zeta_{j,n}-\zeta_{i,n}}$$ $$\leq\frac{c}{n^2}\sum_{j=2}^{n}\sum_{i=1}^{j-1}\frac{1}{j^2-i^2}=O\left(\frac{\log^2 n}{n^{2}}\right).$$

\noindent $M_6=O\left(\frac{1}{n^{2}}\right)$.  
$$M_6=\frac{1}{n(n-1)}\sum_{i=1}^nx_i-\frac{1}{(n-1)(n-2)}\sum_{i=1}^{n-1}\nu_i.$$
Since $r(x)=\prod_{i=1}^n(x-x_i)=(-1)^nn!q_{m,n}=x^n-a_{n-1}x^{n-1}+\dots$, and $r(x)$ fulfils the differential equation $r^{''}+M_{1,n}r^{'}+N_{1,n}r=0$, we can compute
$$a_{n-1}=\sum_{i=1}^nx_i=\lim_{x\to \infty}\frac{x^n-r(x)}{x^{n-1}}=\lim_{x\to \infty}\frac{N_{1,n}(x)x^n+r^{''}(x)+M_{1,n}(x)r^{'}(x)}{N_{1,n}(x)x^{n-1}}.$$
Taking into consideration that $\frac{S^{'}}{S}=\frac{m}{x}$+smaller terms, $\frac{P^{'}}{P}=\frac{m}{x}$+smaller terms and $\frac{P^{''}}{P}=\frac{m(m-1)}{x^2}$+smaller terms, the numerator $N(x)=(n^2-m^2+\alpha(n-m)+(n-1)a_{n-1})x^{n-2}$+smaller terms, and the denominator $D(x)=nx^{n-2}$+smaller terms,
$$a_{n-1}=\lim_{x\to \infty}\frac{N(x)}{D(x)}=\frac{n^2-m^2+\alpha(n-m)+(n-1)a_{n-1}}{n},$$
that is 
$$a_{n-1}=\sum_{i=1}^nx_i=(n-m)(n+m+\alpha).$$
Finally
$$M_6=\frac{(n-m)(n+m+\alpha)}{n(n-1)}-\frac{(n-1-m)(n-1+m+\alpha)}{(n-1)(n-2)}$$ $$=\frac{2m^2+2m\alpha-(\alpha+1)n}{n(n-1)(n-2)}.$$

\noindent $M_1=O\left(\frac{\log n}{n^{2}}\right)$. To compute this let us recall that by \cite{sz} 6.71.12 and 6.71.6
$$\log\left(\prod_{1\leq i<j\leq n}(\zeta_{j,n}-\zeta_{i,n})^2\right)^{\frac{1}{n(n-1)}}=\frac{1}{n(n-1)}\sum_{k=1}^n(k\log k+(k-1)\log(k+\alpha).$$
Let $f(x)=x\log x+(x-1)\log(x+\alpha)$. Since if $x\geq 1$, $f^{'}$ is increasingly tends to infinity at infinity, $0<\frac{f(1)}{2}+\sum_{k=2}^{n-1}f(k)+\frac{f(n)}{2}-\int_1^nf(x)dx<cf^{'}(n)=c\log n$. This entails
$$\sum_{k=1}^n(k\log k+(k-1)\log(k+\alpha)=O(\log n)+\frac{1}{2}n\log n+\frac{1}{2}(n-1)\log (n+\alpha)$$ $$-\frac{n^2}{4}-\frac{(n+\alpha)^2}{4}+\frac{n^2}{2}\log n+\frac{(n-1)^2}{2}\log(n+\alpha)+(1+\alpha)n.$$
So
$$M_1=O\left(\frac{\log n}{n^{2}}\right)+\log\frac{n}{n-1}$$ $$+\frac{1}{4}\left(\frac{n^2+(n+\alpha)^2}{n(n-1)}-\frac{(n-1)^2+(n-1+\alpha)^2}{(n-1)(n-2)}\right)+(1+\alpha)\left(\frac{1}{n-2}-\frac{1}{n-1}\right)$$ $$+\frac{1}{2}\left(\log\frac{n-1+\alpha}{n+\alpha}+\frac{n}{n-2}\log(n-1)-\frac{n+1}{n-1}\log n\right)$$
$$=O\left(\frac{\log n}{n^{2}}\right)+\frac{1}{2}\left(\log\frac{n-1+\alpha}{n+\alpha}+\log\frac{n-1}{n}\right)+\log\frac{n}{n-1}=O\left(\frac{\log n}{n^{2}}\right).$$
So the sum of the error terms is of $O\left(\frac{\log^2 n}{n^{2}}\right)$, and the proof is finished.

\medskip

\section{$X_m$-Laguerre-(II), and $X_m$-Jacobi Polynomials}

The exceptional zeros of exceptional Laguerre-(II) and Jacobi polynomials can be complex, but the asymptotics of these zeros are known. So by the previous methods we can prove asymptotical results and we can examine the regular zeros in these cases.

\subsection{Jacobi Case}

$P_n^{(\alpha,\beta)}(x)$ are the classical Jacobi polynomials, which are orthogonal with respect to the weight $w^{(\alpha,\beta)}(x)=(1-x)^{\alpha}(1+x)^{\beta}$ on the interval $(-1,1)$ ($\alpha, \beta >-1$). $P_n^{(\alpha,\beta)}$ satisfies the differential equation
$$(1-x^2)y^{''}+\left(\beta-\alpha-(\alpha+\beta+2)x\right)y^{'}+ n(n+\alpha+\beta+1)y=0.$$
$\left\{\hat{P}^{(\alpha,\beta)}_{m,m+n}\right\}_{n=0}^{\infty}$ are the exceptional Jacobi polynomials with codimension $m\geq 1$. Let us denote by $S(x):=P_m^{(-\alpha-1,\beta-1)}(x)$. The polynomials $\hat{P}_{m,m+n}^{(\alpha,\beta)}$ are orthogonal with respect to the weight $\frac{w^{(\alpha,\beta)}(x)}{S^2(x)}$ on $(-1,1)$ and satisfy the differential equation
$$y^{''}+\left(\frac{\beta-\alpha-(\alpha+\beta+2)x}{1-x^2}-\frac{2S^{'}(x)}{S(x)}\right)y^{'}$$
\begin{equation}+\left(\frac{m(\alpha-\beta-m+1)+n(n+\alpha+\beta+1)}{1-x^2}-\frac{\beta}{1+x}\frac{2S^{'}(x)}{S(x)}\right)y=0,\end{equation}
(cf. \cite{gumm} (58) (69)). The location of zeros of classical Jacobi polynomials are more complicated then in the Laguerre case, so it is the situation with the exceptional case as well. Here we have

\medskip

\begin{lemma}(\cite{gumm} Prop. 5.3, 5.4, Cor. 5.1.) Let us suppose that $\alpha, \beta$ and $m$ satisfy the condition $\alpha+1-m-\beta \notin \{0,1,\dots ,m-1\}$, and one of the conditions below

\noindent (A) $\beta, \alpha+1-m\in (-1,0)$

\noindent (B) $\beta, \alpha+1-m\in (0,\infty)$. 

Then $\hat{P}_{m,m+n}^{(\alpha,\beta)}$ has exactly $n$ regular zeros (i.e. zeros in $(-1,1)$) which are all simple, and $m$ exceptional zeros (i.e. zeros out of $[-1,1]$); furthermore the regular zeros of $\hat{P}_{m,m+n}^{(\alpha,\beta)}$ approach the zeros of the classical Jacobi polynomials $P_{n}^{(\alpha,\beta)}$, the exceptional zeros of $\hat{P}_{m,m+n}^{(\alpha,\beta)}$ approach the zeros of $P_m^{(-\alpha-1,\beta-1)}$, as $n$ tends to infinity. \end{lemma}

\medskip

With the notation $\hat{w}_{m}^{(\alpha+1,\beta+1)}:=\frac{w^{(\alpha+1,\beta+1)}(x)}{S^2(x)}$, it can be seen that $M_{m,n}(x)=\left(\log \hat{w}_{m}^{(\alpha+1,\beta+1)}\right)^{'}(x)$, and so the first partial derivatives of the logarithm of the energy function with respect to the weight $\hat{w}_{m}^{(\alpha+1,\beta+1)}$ at the zeros of $\hat{P}_{m,m+n}^{(\alpha,\beta)}$ must be zero. By Lemma 1, the diagonal elements of the Hessian can be expressed as
$$-\frac{2}{3}\Phi(\zeta_i)=\frac{1}{6(1-\zeta_i^2)^2}\left(2\left(2\frac{S^{'}}{S}(\zeta_i)(1-\zeta_i^2)+\alpha+\beta+(\alpha-\beta)\zeta_i\right)^2+g(\zeta_i)\right),$$
where
$$g(x)=\left(-(\alpha-\beta)^2-1+(1+2\alpha)(1+2\beta)+\varrho_{m,n}\right)x^2 $$
$$+2(\beta^2-\alpha^2)x-\left((\alpha+\beta)^2+3+(1+2\alpha)(1+2\beta)+\varrho_{m,n}\right),$$
where
$$\varrho_{m,n}=4n(n+\alpha+\beta+1)-8m(m+\beta-\alpha).$$
For $m=1$, $0<\alpha<\beta$ we have $S(x)=\frac{1}{2}\left((\beta-\alpha)x-(\beta+\alpha)\right)$, for which $y_1>1$ (cf. \cite{sz} 4. 21.2), and as it is pointed out in \cite{dy}, $z_1>1$, and $-\Phi(z_1)>0$, $-\Phi(x_i)<0$. Generally the exceptional zeros can be complex, so we concentrate onto the regular zeros.

\medskip

As in the Laguerre case, let $\hat{P}_{m,m+n}^{(\alpha,\beta)}=P_{m,n}q_{m,n}$, where $q:=q_{m,n}$ has the regular, $P:=P_{m,n}$ has the exceptional zeros. Let us define a new weight function again, which depends on $n$ too: let
$$v_{m,n}^{(\alpha+1,\beta+1)}=P_{m,n}^2\hat{w}_{m}^{(\alpha+1,\beta+1)}.$$ Let us consider the energy function $T_v(u_1,\dots ,u_n)$ on $(-1,1)^n$ with respect to $v:=v_{m,n}^{(\alpha+1,\beta+1)}$. With these notations we have

\medskip

\begin{theorem} If $\alpha>m-1$, $\beta>0$ and if $n$ is large enough, then the set of the regular zeros of $\hat{P}_{m,m+n}^{(\alpha,\beta)}$ is the unique set of minimal energy with respect to the external field $\left(v_{m,n}^{(\alpha+1,\beta+1)}\right)^{\frac{1}{2(n-1)}}$.\end{theorem}

\medskip

\proof Repeting the first part of the proof of Theorem 2, we can see that the first partials of the logarithm of the energy function are zero at the regular zeros. The non-diagonal elements of the Hessian are the same as well. Denoting the zeros of $S,P,q$ by $y_1,..y_m$, $z_1,..z_m$ , $x_1,..x_n$ respectively ($z_i$ are not necessarily distinct and real), as in Th. 2,
$$H_{i,i}=2\left(\frac{P^{'}}{P}-\frac{S^{'}}{S}\right)^{'}(u_i)+\left(\log w ^{(\alpha+1,\beta+1)}\right)^{''}(u_i)-2\sum_{1\leq j\leq n \atop j\neq i}\frac{1}{(u_i-u_j)^2} $$
$$=2\sum_{j=1}^m\left(\frac{1}{(u_i-y_j)^2}-\frac{1}{(u_i-z_j)^2}\right)$$ $$-\frac{\alpha+\beta+2}{(1-u_i^2)^2}\left(u_i^2+\frac{2(\alpha-\beta)}{\alpha+\beta+2}u_i+1\right)-2\sum_{1\leq j\leq n \atop j\neq i}\frac{1}{(u_i-u_j)^2}. $$
Because $\frac{\alpha+\beta+2}{(1-x^2)^2}\left(x^2+\frac{2(\alpha-\beta)}{\alpha+\beta+2}x+1\right)\geq \frac{4(\alpha+1)^2}{\alpha+\beta+2}$ and both in the cases (A) and (B) $\alpha+\beta+2>0$, the second and third terms are negative and uniformly bounded for arbitrary sets of nodes. If $n$ is large enough, by Lemma 4 we can estimate the first sum as
 $$2\sum_{j=1}^m\frac{\left|z_j-y_j\right|\left(\left|u_i-z_j\right|+\left|u_i-y_j\right|\right)}{\left|u_i-z_j\right|^2\left|u_i-y_j\right|^2}\leq \frac{4m\varepsilon_nc_1(m)}{c_2(m)^4}.$$
Here $\left|z_j-y_j\right|\leq \varepsilon_n$, which is small when $n$ is large enough (cf. Lemma 4.), $c_1(m):=2\max_jd(y_j,[-1,1])$, $c_2(m):=\frac{1}{2}\min_jd(y_j,[-1,1])$, where $d$ means distance. The computations above entails that if $n$ is large enough, then $H_{i,i}<0,$ and we can finish the proof as in Th. 3.

\subsection{Laguerre II Case}

 The type II exceptional Laguerre polynomials (cf. \cite{gumm}) of codimension $m\geq 1$ are $\{L_{m,m+n}^{II, (\alpha)}\}_{n=0}^{\infty}$. Let us denote by $S(x):=L_m^{(-\alpha-1)}(x)$. Let us assume that $\alpha>m-1$. By the assumption (cf. \cite{gumm} Prop.4) $S$ has no zeros in $[0,\infty)$ and $L_{m,m+n}^{II, (\alpha)}$-s are orthogonal on $(0,\infty)$ with respect to the weight $\hat{w}^{(\alpha)}:=\frac{x^{\alpha}e^{-x}}{s^2(x)}$ (\cite{gumm} (50)). $L_{m,m+n}^{II, (\alpha)}$ satisfies the following differential equation (\cite{gumm} (40) (31)):
 \begin{equation}y^{''}(x)+\left(\frac{\alpha+1-x}{x}-\frac{2S^{'}(x)}{S(x)}\right)y^{'}(x)+\left(\frac{n-m}{x}-\frac{\alpha}{x}\frac{2S{'}(x)}{S(x)}\right)y(x)=0.\end{equation}
 
 \medskip
 
 \begin{lemma}(\cite{gumm} Prop. 4.3, 4.4, 4.5, 4.8,  Cor. 4.1) $L_{m,m+n}^{II, (\alpha)}$ has $n+m$ simple zeros, $n$ regular zeros  $x_{m,1}^{(\alpha)}, \dots ,x_{m,n}^{(\alpha)} \in (0,\infty)$ and $0$ or $1$ negative zero according to whether $m$ is even or odd. Furtheremore
 $$\lim_{n\to \infty}nx_{m,i}^{(\alpha)}=\frac{\left(j^{(\alpha)}_i\right)^2}{4},$$
 and as $n \to \infty$ the exceptionl zeros of $L_{m,m+n}^{II, (\alpha)}$ converge to the zeros of $S(x)=L_m^{(-\alpha-1)}(x)$.
 \end{lemma}
 
\medskip

We can write $L_{m,m+n}^{II, (\alpha)}=P(x)q(x)$ again, where $P(x)$ is a polynomial of degree $m$ and with $0$ or $1$ real zero. So as in the Jacobi case we can examine the properties of the regular zeros of $L_{m,m+n}^{II, (\alpha)}$, ie the zeros of $q$. As previously, the following if valid:
$$ q^{''}(x)+\left(M_n(x)+2\frac{P^{'}}{P}(x)\right)q^{'}(x)+\left(N_n(x)+\frac{P^{''}}{P}(x)+M_n(x)\frac{P^{'}}{P}(x)\right)q(x)=0,$$
where $M_n(x)$ and $N_n(x)$ are the coefficients in the equation (12). Let $v(x)=v^{(\alpha+1)}_{m,n}(x)=\frac{x^{\alpha}e^{-x}P^2(x)}{s^2(x)}$, as above. With the notations above we have

\medskip

\begin{theorem} If $\alpha>m-1$, and if $n$ is large enough, then the set of the regular zeros of $\hat{L}{m,m+n}^{(II,\alpha)}$ is the unique set of minimal energy with respect to the external field $\left(v_{m,n}^{(\alpha+1)}\right)^{\frac{1}{2(n-1)}}$.\end{theorem}

\medskip

\proof
As in Th.2 and Th.5 the first partial derivatives of $\log T_v$ are zero at the regular zeros of $L_{m,m+n}^{II, (\alpha)}$, and with the previous notations for any set $U_n\subset(0,\infty)$
$$H_{i,i}=2\left(\frac{P^{'}}{P}-\frac{S^{'}}{S}\right)^{'}(u_i)-\frac{\alpha+1}{u_i^2}-2\sum_{1\leq j \leq n\atop j\neq i}\frac{1}{(u_i-u_j)^2},$$
where by Lemma 5 the first term tends to zero as $n\to \infty$ and the last two terms are negative. Finally we can finish the proof as in Th. 2 again.

\medskip

\noindent \small{Department of Analysis, \newline
Budapest University of Technology and Economics}\newline
\small{ g.horvath.agota@renyi.mta.hu}

\end{document}